\documentclass[12pt]{article}
\usepackage{theorem,amsfonts,amssymb}
\newtheorem{theorem}{Theorem}[section]
\newtheorem{proposition}{Proposition}[section]
\newtheorem{lemma}{Lemma}[section]
\newtheorem{corollary}{Corollary}[section]
\theorembodyfont{\upshape}
\newtheorem{definition}{Definition}
\newtheorem{remark}{Remark}

\newtheorem{example}{Example}[section]
\newtheorem{proof}{Proof}

\newtheorem{acknowledgement}{Acknowledgement}

\newcommand{\bt}{\begin{theorem}}
\newcommand{\et}{\end{theorem}}
\newcommand{\bl}{\begin{lemma}}
\newcommand{\el}{\end{lemma}}
\newcommand{\bp}{\begin{proposition}}
\newcommand{\ep}{\end{proposition}}
\newcommand{\bex}{\begin{example}}
\newcommand{\eex}{\end{example}}
\newcommand{\bc}{\begin{corollary}}
\newcommand{\ec}{\end{corollary}}
\newcommand{\bo}{\begin{proof}}
\newcommand{\eo}{\end{proof}}
\newcommand{\bd}{\begin{definition}}
\newcommand{\ed}{\end{definition}}
\newcommand{\br}{\begin{remark}}
\newcommand{\er}{\end{remark}}
\newcommand{\be}{\begin{enumerate}}
\newcommand{\ee}{\end{enumerate}}

\newcommand{\R}{{\mathbb R}}
\newcommand{\Z}{{\mathbb Z}}
\newcommand{\Q}{{\mathbb Q}}
\newcommand{\N}{{\mathbb N}}
\newcommand{\cI}{{\cal I}}

\begin{document}

\title{Dynamic Random Walks on Motion Groups}
\author{C. R. E. Raja and R. Schott}
\date{}
\maketitle

\let\epsi=\epsilon
\let\vepsi=\varepsilon
\let\lam=\lambda
\let\Lam=\Lambda
\let\ap=\alpha
\let\vp=\varphi
\let\ra=\rightarrow
\let\Ra=\Rightarrow
\let \Llra=\Longleftrightarrow
\let\Lla=\Longleftarrow
\let\lra=\longrightarrow
\let\Lra=\Longrightarrow
\let\ba=\beta
\let\ov=\overline
\let\ga=\gamma
\let\Ba=\Delta
\let\Ga=\Gamma
\let\Da=\Delta
\let\Oa=\Omega
\let\ol=\overline
\let\Lam=\Lambda
\let\un=\upsilon

\begin{abstract}
In this note, we give an original convergence result for products of
independent random elements of motion group. Then we consider
dynamic random walks which are inhomogeneous Markov chains whose
transition probability of each step is, in some sense, time
dependent. We show, briefly, how Central Limit theorem and Local
Limit theorems can be derived from the classical case and provide
new results when the rotations are mutually commuting. To the best
of our knowledge, this work represents the first investigation of
dynamic random walks on the motion group.
\end{abstract}



\section{Introduction}
The motion group $G=SO(d)\ltimes \R ^d$ is the semi-direct product
of $SO(d)$, the group of rotations in the space $\R ^d$, and $\R
^d$. This group plays a special role in the study of random walks on
Lie groups \cite{GKR}. A Central Limit Theorem on motion groups has
been proved by Roynette \cite{ro} and Baldi, Bougerol, Crepel
\cite{bbc} gave a Local Limit Theorem on Motion groups. Random walks
on homogeneous spaces of the motion group have been studied by
Gallardo and Ries \cite{gr}. The main novelty of this paper is in
the dynamic model of random walks which we define on the motion
group. The theory of dynamic random walks has been done by
Guillotin-Plantard in a commutative setting \cite{gu1, gu2, gu3}. So
far, dynamic random walks have been considered on Heisenberg groups,
the dual of $SU(2)$ \cite{GPS} and Clifford algebras \cite{sss}.
Needless to say, there is much work to do. This paper is another
(small) step/attempt in extending this theory to non-commutative
algebraic structures. Recently, random walks on the motion group
have been proposed as algorithmic ingredients for searching in
peer-to-peer wireless networks \cite{gms}. The organization is as
follows: Section II contains basic definitions, known limit
theorems, as well as a new convergence theorem for product of random
elements of the motion group. Dynamic random walks are considered in
Section III, we recall known results, show how to derive some limit
theorems from the classical case and investigate more deeply when
the rotations form an abelian group.  It may be noted that these are
the first examples of non-discrete dynamic random walks.  Section IV
provides some concluding remarks and further research aspects.

\section{Motion group}
\subsection{Basic definitions and known results}

The composition law of the motion
group $G=SO(d)\ltimes \R ^d $ is given by:
$$(R_1,T_1).(R_2,T_2)=(R_1\circ R_2, T_1+R_1(T_2))$$
Remember that $R_1\circ R_2$ is the rotation of angle $\Theta_1+\Theta_2$
if $\Theta_1$ (resp. $\Theta_2$) is the rotation angle of $R_1$
(resp.$R_2$).
More generally:$$(R_1,T_1)(R_2,T_2)...(R_n,T_n)=(R_1\circ
R_2\ldots\circ R_n,
T_1+R_1(T_2)+R_1\circ R_2(T_3)+...+R_1\circ R_2\ldots\circ R_{n-1}(T_n))$$
where $(R_i , T_i)$ are $G$-valued random variables.
Let $$S_n=T_1+R_1(T_2)+R_1\circ R_2(T_3)+...+R_1\circ R_2\ldots\circ R_{n-1}(T_n)$$
$S_n$ gives the position after n steps. $S_n$ is the sum of $n$
(not necessarily independent) random variables.\\
The following Central Limit Theorem has been proven in \cite{ro}:
\begin{theorem}
Assume that $R_i$, (resp. $T_i$), $i\in\{1,2,\ldots,n, \ldots \}$
are $n$ independent random variables with common law $\mu$ (resp
$\nu$), that the support of  $\mu$ generates $SO(d)$ and that $\nu$
has a second
  order moment. Then $\frac{S_n}{ \sqrt{n}}$ converges in law to
  the Gaussian distribution $N(0,{\theta I_d})$ when $n$ goes to
  infinity. $I_d$ stands for the $d\times d$ dimensional identity matrix and
  $\theta$ is a positive constant.
\end{theorem}

{\bf Remark:} This theorem tells us, intuitively, that $\frac{S_n}{ \sqrt{n}}$
becomes rotation invariant when $n$ goes to infinity and that $S_n$
behaves asymptotically as a random walk on ${\R}^d$ which is rotation
invariant. In other words:
$$S_n\sim_{n\rightarrow\infty} Y_1+Y_2+\ldots +Y_n$$
where $Y_i$, $i\in\{1,2,\ldots,n\}$ are $n$ independent and
identically distributed random variables.

The following Local Limit Theorem has been proven in \cite{bbc}, we
formulate it below in a simple way:
\begin{theorem}
Let $P_n(O,M)$ be the probability that the random walks on $G$
reaches the point M of $\R ^d$ in n steps when staring from the
point O, then::
$$P_n(O,M)=P(S_n=M)\sim_{n\rightarrow\infty}\frac{K}{n^{d/2}}$$
where $K$ is a positive constant (independent of n).
\end{theorem}
\subsection{A convergence theorem}
Let $O(d)$ be the group of orthogonal linear transformations on $\R
^d$ ($d\geq 1$) and $K$ be a compact subgroup of $O(d)$ and $G = K
\ltimes \R ^d$ be a group of motions of $\R ^d$.  Let $Y_i = (R_i,
T_i)$ be independent random variables. Let $S_n = T_1+R_1T_2 +
\ldots + R_1 R_2 \ldots R_{n-1} T_n$ and $X_n = R_1 R_2 \ldots
R_{n-1} T_n$ with $R_0 = 1$.

\bt\label{lln}
Assume the following:
\be
\item $R_1 \ldots R_n$ converges in law to $\omega _K$, the normalized
Haar measure on $K$;

\item $\R ^d$ has no non-zero $K$-invariant vectors;

\item $X_n$ has second moment;

\item $E(T_n)$ is bounded.
\ee
Then
$${S_n\over b_n}\to 0~~ a.s$$
for any sequence $(b_n)$ such that $\sum {E(||X_n-E(X_n)||^2)\over
b_n ^2}<\infty$. \et

\bo We recall that a random vector $T$ in $R^d$ is said to have
finite expectation if there is a vector $v\in \R^d$ such that $<v,u>
= E(<T, u>)$ for any $u \in \R ^d$ and in this case we write $E(T) =
v$.  Also if $R$ is a random rotation on $\R ^d$, then $E(R)$ is a
operator on $\R ^d$ defined by
$$<E(R)u, v> = E(<Ru, v>)$$ for any two vectors $u, v \in \R ^d$.

It follows that $E(X_n) = E( R_1 R_2 \ldots R_{n-1}) E(T_n)$ for all
$n \geq 1$. For $u , v\in \R^d$,
$$<E(R_1R_2\ldots R_n)u, v> = \int <R_1R_2\ldots R_n u, v> d\omega =
\int <T(u), v> \rho _n (dT)$$ where $\rho _n$ is the distribution of
$R_1R_2\ldots R_n$. Since $R_1 R_2 \ldots R_n $ converges in law to
$\omega _K$, we get that $E(R_1 R_2 \ldots R_{n-1}) \to P_K$ where
$P_K$ is the projection onto $K$-fixed vectors.

We first claim that $E(X_n) \to 0$.  Since $\R ^d$ has no
$K$-invariant vectors, $P_K = 0$. Now since $E(T_n)$ is bounded, if
$v$ is a limit point of $E(T_n)$, let $E(T_{k_n}) \to v$.  Then
since $E(R_1 R_2 \ldots R_{n-1}) \to 0$ in the operator topology,
$E(X_{k_n}) \to 0$. Thus, $E(X_n) \to 0$.

Let $u \in \R ^d$.  Take $Z_n = <X_n -E(X_n), u>$.  Then $E(Z_n )
=0$. $Z_n$ are independent real random variables with finite second
moments. Then $${1\over b_n} \sum _{i=1} ^n Z_i \to 0~~ a.s.$$ for
any constant $(b_n)$ such that $\sum _0 ^\infty {Var(Z_n) \over
b_n^2} <\infty$ (cf. \cite{S}).  This implies that
$${1\over b_n} \sum _{i=1}^n(X_i - E(X_i)) \to 0 ~~a.s.$$   We have shown that
$E(X_n) \to 0$ and hence ${1\over b_n}\sum _{i=1}^n E(X_i) \to 0$.
Thus,
$${1\over b_n}\sum _{i=1}^n X_i \ra 0 ~~a.s.$$
\eo

The conditions in Theorem \ref{lln} are verified if we take $R_i$ to
be iid with the support of the common law is aperiodic (that is,
support is not contained in a coset of a proper normal subgroup) and
$T_i$ to be dynamic random walk with $b_n= {1\over n^\ap}$ for any
$\ap >{1\over 2}$. Thus, under these assumptions we get that
$${1\over n^\ap}(T_1+R_1T_2 + \ldots + R_1 R_2 \ldots R_{n-1} T_n)
\to 0 ~~a.s.$$

\section{Dynamic random walks}
\subsection{Preliminaries and known results}
Let $S=(E,{\cal A},\rho ,\tau )$ be a dynamical system
\index{dynamical
  system} where $(E,{\cal A},\rho)$ is a probability space \index{probability space} and
$\tau $ is a transformation defined on $E$. Let $d\geq 1$ and
$h_{1},\ldots,h_{d}$ be functions defined on $E$ with values in
$[0,\frac{1}{d}]$. Let $(T_{i})_{i\geq 1}$ be a sequence of
independent random vectors with values in $\Z ^{d}$. Let $x\in E$
and $(e_{j})_{1\leq j\leq d}$ be the unit coordinate vectors of $\Z
^{d}$. For every $i\geq 1$, the law of the random vector $T_{i}$ is
given by
$$P(T_{i}=z)=\left\{\begin{array}{ll}

 h_{j}(\tau ^{i}x) & \mbox{if}\ \ z=e_{j}\\

\frac{1}{d}-h_{j}(\tau ^{i}x) & \mbox{if}\ \ z=-e_{j}\\

0 & \mbox{otherwise}

\end{array} \right.
$$
We write $$S_{0}=0,\ \ S_{n}=\sum_{i=1}^{n} T_{i}\ \mbox{for}\ n\geq
1$$ for the ${\Z}^{d}$-random walk generated by the family
$(T_{i})_{i\geq 1}$. The random sequence $(S_n)_{n\geq 0}$ is called
a dynamic
 ${\Z}^{d}$-random walk.

It is worth remarking that if the functions $h_{j}$ are constant
then we have the classical random walks but if these functions are
all not constant, $(S_{n})_{n\in \N}$ is a non-homogeneous
Markov chain.\\

Let ${\cal C}_{1}(S)$ denote
the class of functions $f\in L^{1}(E,\mu)$ satisfying the
following condition $(H_1)$:
$$\left|\sum_{i=1}^{n}\Big(f(\tau ^{i}x)-\int_{E}f(x) d\rho (x)\Big)\right|=
o\Big(\frac{\sqrt{n}}{log(n)}\Big)$$

Let ${\cal C}_{2}(S)$ denote
the class of functions $f\in L^{1}(E,\mu)$ satisfying the
following condition $(H_2)$:
$$\sup_{x\in E}\left|\sum_{i=1}^{n}\Big(f(\tau ^{i}x)-\int_{E}f(x) d\rho (x)\Big)\right|=
o\Big(\sqrt{n}\Big)$$

A Central Limit Theorem:
\begin{theorem}\label{tcl}
Asssume that for every $j,l\in\{1,\ldots,d\}$, $h_{j}\in{\cal
C}_2(S)$, $h_{j}h_{l}\in{\cal C}_2(S)$ and
$\int_{E}h_{j}d\rho=\frac{1}{2d}$. Then, for every $x\in E$, the
sequence of processes $(\frac{1}{\sqrt{n}}S_{[nt]})_{t\geq 0}$
 weakly converges in the Skorohod space ${\cal D}={\cal D}([0,\infty[)$  to the $d$-dimensional
Brownian motion
$$B_{t}=(B_{t}^{(1)},\ldots,B_{t}^{(d)})$$
with zero mean and covariance matrix $A t$.\index{Brownian motion}
\end{theorem}
The proof of this theorem is in \cite{GPS}.\\

A Local Limit Theorem:
\begin{theorem}
Let $h_{j}\in{\cal C}_1(S)$, $h_{j}h_{l}\in{\cal C}_1(S)$ and
$\int_{E}h_{j}d\rho =\frac{1}{2d}$. Then, for almost every $x\in E$,
$P(S_{2n}=0)$, the probability that starting from the point O, the
random walks comes back to O in $2n$ steps, has the following
asymptotic behavior:
$$P(S_{2n}=0)\sim\frac{2}{\sqrt{det(A)}(4\pi n)^{d/2}}$$
as $n\rightarrow\infty$

\end{theorem}

The proof of this theorem is also in \cite{GPS}.\\

\subsection{Dynamic random walks on the motion group}
Recall that we consider the random walk
$$S_n=T_1+R_1(T_2)+R_1\circ R_2(T_3)+...+R_1\circ R_2\ldots\circ R_{n-1}(T_n)$$
where $T_i$, $i\in{\N}$ are dynamic random variables as
defined above and we now define dynamic random rotations $R_i$.\\

If the rotations are classical random variables and translations are
dynamic random variables then one can adapt the result in \cite{ro}
and prove a Central Limit Theorem and a Local Limit Theorem
\cite{bbc} for $S_n$ thanks to the Central Limit Theorem and the
Local Limit Theorem for dynamic random walks \cite{GPS} given in the
above section. We do not write explicitely these theorems because
these formulation is almost the same as in \cite{bbc}, \cite{ro}.
Similar Central Limit Theorem and Local Limit Theorem hold true
under Lindenberg
condition on the translations $T_i$.\\

If both rotations and translations are dynamic random walks the
problem is still open. \\

We consider now the $2$-dimensional case. It is a known that $SO(2)$
is a compact abelian group (isomorphic to $U(1)$) and for any
irrational number $\theta \in \R$, $e^{2\pi i\theta }$ generates a
dense subgroup of $SO(2)$.  Using this fact we prove that the
convolution product $\mu_1*\mu_2*\ldots*\mu_n$ of dynamic measures
corresponding to dynamic rotations $R_1$,\ldots, $R_n$ converges
weakly to the Haar measure of $SO(2)$.

Let $\theta$ be an irrational number and $R_j$ be random rotations
on $\R ^2$ defined by
$$P(R_{j}=z)=\left\{\begin{array}{ll}

 f(\tau ^{j}x) & \mbox{if}\ \ z=e^{2\pi i\theta }\\

1-f(\tau ^{j}x) & \mbox{if}\ \ z=e^{-2\pi i \theta}\\

0 & \mbox{otherwise}

\end{array} \right.$$
and $f\colon E \to [0, 1]$ satisfies $f(1-f)\in {\cal C}_2(S)$ where
$E$ and ${\cal C}_2(S)$ are as in 3.1 with $d=1$.

If $f$ is an indicator function taking values $0$ and $1$, then it
can be easily seen that $R_i$ degenerate and hence the product
$R_1\ldots R_n$ does not converge (in law) as the set $\{ e^{2\pi i
k \theta} \mid k\geq 1 \}$ is dense in $SO(2)$.  This forces us to
assume that $f$ is not a indicator function.   In this case, we have
the following:

\bt\label{dr} Almost surely $R_1R_2 \ldots R_n $ converges in law to
the Haar measure on $SO(2)$. \et

In order to prove the above result we need to recall some details on
the dual of compact abelian groups and Fourier transform of
probability measures on compact groups.

\noindent {\bf Dual of compact groups:} For a compact abelian group
$K$, continuous homomorphisms from $K$ into $SO(2)$ are known as
characters and characters form a (locally compact abelian) group
which is denoted by $\hat K$ and is called the dual group of $K$:
cf. \cite{HR} for details on duality of locally compact abelian
groups.  For each integer $n$, the map $z\mapsto z^n$ defines a
character on $SO(2)$ and defines an isomorphism between the group
$\Z$ of integers with the dual of $SO(2)$ (cf. 23.27 (a) of
\cite{HR}).  It is known that if $K_1, \ldots , K_d$ are compact
abelian groups, then the dual of $K_1\times \ldots \times K_d$ is
isomorphic to $\hat K_1 \times \ldots \times \hat K_d$ (cf. 23.18 of
\cite{HR}).

\noindent {\bf Fourier transform:} Let $K$ be a compact abelian group and
$\mu$ be a probability measure on $K$.  Then the Fourier transform of
$\mu$, denoted by $\hat \mu$ is a function on $\hat K$ and is defined by
$$\hat \mu (\chi ) = \int \chi (x) d\mu (x)$$ for all $\chi \in \hat K$.
It is known that $\mu $ is the normalized Haar measure on $K$ if and only
if
$$\hat \mu (\chi )=\left\{\begin{array}{ll}

 0 & \mbox{if}\ \ \chi {\rm ~~ is ~~non-trivial}\\

1  & \mbox{if}\ \ \chi {\rm ~~ is ~~trivial}

\end{array} \right.
$$
and if $X_n$ are $K$-valued random variables with Fourier transform
$f_n$, then $X_n$ converges in law to a $K$-valued random variable
$X$ if and only if $f_n $ converges to the Fourier transform of $X$
pointwise (cf. \cite{Ru}).

\bo $\!\!\!\!\!$ {\bf of Theorem \ref{dr}}\ \ Let $k$ be any
non-zero integer.  It is sufficient to claim that
$$\prod _{j=1}^n|\int e^{2\pi i kx }d\mu _j |\to 0$$ as $n \ra
\infty$.

$$\begin{array}{lcl}
|\int e^{2\pi i kx }d\mu _j |^2 & =&|e^{2\pi ik\theta } f(\tau ^jx)
+ e^{-2\pi i k\theta}(1-f(\tau ^jx)) |^2 \\ &=& |\cos (2\pi k\theta)
+i \sin (2\pi k\theta )(f(\tau ^jx)-1+f(\tau ^jx)) |^2 \\
& = & \cos ^2(2\pi k\theta)+\sin ^2(2\pi k\theta )(1-2f(\tau ^jx)) ^2\\
& = & 1- 4\sin ^2 (2\pi k\theta )f(\tau ^jx)(1-f(\tau ^jx))\\
\end{array}$$

Suppose $f(\tau ^jx) (1-f(\tau ^jx) )\not \to 0$.  Then $1- 4\sin ^2
(2\pi k\theta )f(\tau ^jx)(1-f(\tau ^jx))\not \to 1$ and hence
$\prod _{j=1} ^n |\int e^{2\pi i k\theta } d\mu _j | \to 0$.   Thus,
it is sufficient to show that $f(\tau ^jx)(1-f(\tau ^jx))\not \to
0$.

If $f(\tau ^jx) (1-f(\tau ^jx)) \to 0$, then
$${1\over n}\sum _1 ^n f(\tau ^jx) (1-f(\tau ^jx)) \to 0 = \int
f(x) (1-f(x))d\rho (x)$$ and hence $f$ is an indicator function.
This is a contradiction.  Thus proving the result. \eo

Let $K$ be a compact connected abelian subgroup of $SO(d)$, for
instance one may take $K$ to be the maximal torus in $SO(d)$.  In
this situation one can define dynamic random walks in many ways and
we will now consider two forms of dynamic random walks on $K$. The
first one is the following: $a \in K$ is such that the closed
subgroup generated by $a$ is $K$ (see 25.15, \cite{HR} for existence
of such $a$) and $R_j$ are random variables taking values in $K$
defined by
$$P(R_{j}=x)=\left\{\begin{array}{ll}

 f(\tau ^{j}x) & \mbox{if}\ \ x=a\\

1-f(\tau ^{j}x) & \mbox{if}\ \ x=a^{-1}\\

0 & \mbox{otherwise}

\end{array} \right.$$
and $f\colon E \to [0, 1]$ satisfies $f(1-f)\in {\cal C}_2(S)$ where
$E$ and ${\cal C}_2(S)$ are as in Section 3.

In the situation we have the following as a consequence of Theorem
\ref{dr}.

\bt\label{ac1} Almost surely $R_1R_2 \ldots R_n $ converges in law
to the Haar measure on $K$. \et

\bo For any non-trivial character $\chi $ on $K$, the element $\chi
(a)$ in $SO(2)$ corresponds to an irrational rotation, hence we get
from Theorem \ref{dr} that $\chi (R_1 R_2 \ldots R_n)$ converges in
law to the Haar measure on $SO(2)$ which is $\chi (\omega _K)$. This
implies that $(R_1 R_2 \ldots R_n)$ converges in law to the Haar
measure on $K$ \eo

\br The Corollary \ref{ac1} could be proved for any monothetic
compact connected abelian group in a similar way but for simplicity
and for the purpose of the article we restrict our attention to
compact connected abelian subgroups of $SO(d)$: a topological group
$K$ is called monothetic if $K$ contains an element $a$ such that
the closed subgroup generated by $a$ is $K$ itself (cf. 9.2 of
\cite{HR} for monothetic compact groups). \er

We will now consider the second form of dynamic random walks on $K$.
Let $v_1,\cdots , v_r$ be a basis for the Lie algebra of $K$ and
$\exp$ be the exponential map of the Lie algebra of $K$ into $K$.
Let $e_k = \exp (v_k)$ $1\leq k \leq r$.  Let $R_j$ be the random
variables taking values in $K$ defined by

$$P(R_{j}=x)=\left \{ \begin{array}{ll}
f_k(\tau ^{j}x) & \mbox{if}\ \ x=e_k\\
\frac{1}{r}-f_k(\tau ^{j}x) & \mbox{if}\ \ x=e_k^{-1}\\
0 & \mbox{otherwise}
\end{array} \right.$$
and $f_k$ are functions from $E$ taking values in $[0, {1\over r}]$
where $E$ is as in Section 3.  We further assume that $k$-the
coordinate of $v_k$ is irrational so that $e_k$ is an irrational
rotation by an angle $\theta _k$ and all other coordinates of $v_k$
are $0$. We further assume that $1$ and $\theta _k$ are independent
over $\Q$.

In this situation also we have the following which could be proved
as in Theorem \ref{ac1}

\bt\label{ar} Almost surely $R_1R_2 \ldots R_n $ converges in law to
the Haar measure on $K$. \et

As an application of the results proved in the previous section and
the above results on compact groups we get the following:

\bt\label{ct} Let $(R_j, T_j)$ be dynamic random walk on $K\ltimes
\R ^{d}$ where $R_j$ is the dynamic random walk on $K$ given in
Theorem \ref{ac1} or Theorem \ref{ar} and $T_j$ is dynamic random
walk on $\R ^{d}$ defined in 3.1.  Then for $\ap >{1\over 2}$,
$${1\over n^\ap}(T_1+R_1T_2 + \ldots + R_1 R_2 \ldots R_{n-1} T_n)
\to P_K(v _0) ~~a.s$$ where $P_K$ is the projection onto the
$K$-invariant vectors in $R^d$ and $v_0= (2E(h_j|\cI)-{1\over
d})_{1\leq j \leq d}$.
\et

\bo We first assume that $\R ^d$ has no non-zero $K$-invariant
vectors.  Condition (1) of Theorem \ref{lln} follows from Theorems
\ref{ac1} and \ref{ar}. Let $X_n = R_1 R _2 \ldots R_{n-1} T_n$.
Then $E(<X_n, u>^2) = \int <R_1R_2\ldots R_{n-1}T_n, u> ^2 d\omega $
is finite as $T_n$ takes only finitely many values and rotations
preserve the norm. Thus, Condition (3) of Theorem \ref{lln} is
verified and condition (4) is easy to verify.  Hence $${1\over
n^\ap}(T_1+R_1T_2 + \ldots + R_1 R_2 \ldots R_{n-1} T_n) \to 0
~~a.s$$

In general, let $V$ be the space of $K$-invariant vectors in $\R
^d$.  Let $P_K$ be the orthogonal projection onto $V$ and $Q$ be the
projection onto the orthogonal complement of $V$. Then for any $v\in
\R ^d$, $v = P_K(v) +Q(v)$ and $Q(\R ^d)$ has no non-zero
$K$-invariant vector.  Since both $V$ and $Q(\R^d)$ are
$K$-invariant, we get that $P_K(R(v))=P_K(v)$ $Q(R(v)) = R(Q(v))$
for any $v \in \R ^d$ and $R \in K$.  Now the result follows from
the previous case and Theorem 2.1 of \cite{GPS}

\eo

\section{Concluding remarks}
We have proved a new convergence result for classical random walks on
the motion group. Our results for the dynamic case are still partial
and we are planning to characterize recurrent and transient random
walks (in this model) on the motion group and the corresponding
homegeneous spaces.
So far, dynamic random walks have only been considered on Heisenberg
groups, the dual of $SU(2)$ \cite{GPS}, the motion group and Clifford
algebras \cite{sss}. A more general study of dynamic random walks
on Lie groups, homogeneous spaces and quantum groups is still to be
done. This is a challenging research project.

\begin{acknowledgement}
This work is done as a part of the IFIM project on "Dynamic random
walks on algebraic structures and applications in quantum
computing".  The authors would like to thank IFIM for its support.
\end{acknowledgement}

\noindent{C. R. E. Raja, Stat-Math Unit, 8th Mile Mysore Road,
Bangalore 560 009, India.}

\noindent{e-mail: creraja@isibang.ac.in}

\medskip

\noindent{R. Schott, IECN and Loria, Univesit\'e Henri Poincar\'e,
BP 239, Nancy, Rrance.}

\noindent{e-mail: schott@loria.fr}

\end{document}